\documentclass[12pt,a4paper,oneside]{article}
\usepackage[T1]{fontenc}
\usepackage[utf8]{inputenc}
\usepackage{amsmath,amssymb}
\usepackage[russian]{babel}
\usepackage{graphicx}
\usepackage{float}
\usepackage{mathtools}
\usepackage{biblatex}
\usepackage[a4paper,left=2.5cm,right=2.5cm,top=2.5cm,bottom=2.5cm]{geometry}
\bibliography{pi-blocks}
\title{\textbf{Получение $\pi$ через столкновение брусков}}
\author{Иван Людвиг Терешко\\МФТИ}
\date{}
\begin{document}
\maketitle
\begin{abstract}
В работе рассматривается метод получения числа $\pi$ как количество упругих столкновений между двумя брусками и стеной.
\end{abstract}

\section{Введение}
\begin{figure}[H]
	\includegraphics[width=\linewidth]{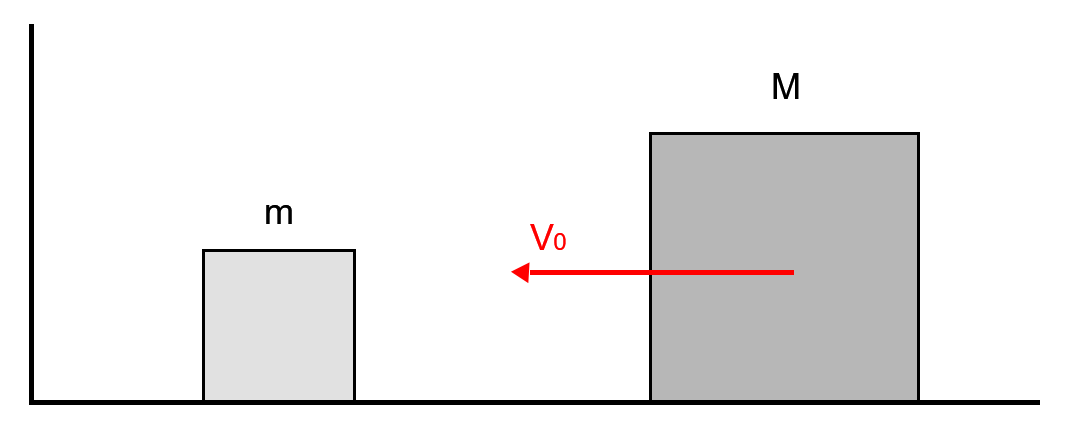}
	\caption{Рассматриваемая модель}
	\label{model}
\end{figure}

Число $\pi$ встречается в некоторых физических законах. В данной статье рассматривается простая модель, в которой оно обнаруживается. Модель состоит из гладкой поверхности, двух брусков, стены (Рис. \ref{model}). Все столкновения абсолютно упругие. Стол трения между брусками и поверхностью нет. Малый брусок имеет массу $m$ расположен между стенкой и большим бруском. В начальный момент времени больший брусок массой $M$ имеет скорость $V_0$, направленную в сторону стены. Далее, происходят столкновения между брусками, а также со стенкой.\\
\begin{table}[h]
	\centering
	\begin{tabular}{ | c  | c | }
		\hline
		$\alpha$ & $N$ \\ \hline
		$1$ & 3 \\ \hline
		$10^{-2}$ & 31 \\ \hline
		$10^{-4}$ & 314 \\ \hline
		$10^{-6}$ & 3141 \\ \hline
		$10^{-12}$ & 3141592 \\ \hline
	\end{tabular}
	\label{table}
	\caption{Количество столкновений}
\end{table}
Рассмотрим количество столкновений в системе при различных отношениях масс брусков $\alpha = m/M$ \cite{3b1b}. Из данных видно, что при отношении масс $10^{-2n}$ количество столкновений равно числу, состоящему из первым $n$ цифр числа $\pi$.\\
Сандерсон (2019) опубликовал видео \cite{3b1b}, посвящённое данной задаче, которую впервые рассмотрел Галперин (2003). В этой работе я использую другой подход для решения задачи и сравниваю его с ранее предложенным. 

\section{Описание}
Рассмотрим столкновения между брусками. $m$, $M$, $v$, $V$ — массы малого и большего бруска и их скорости до столкновения соответственно. $v'$ и $V'$ — скорости малого и большего бруска после столкновения соответственно.
Удар абсолютно упругий, значит выполняется закон сохранения энергии:
\begin{align}
	\label{coe}
	\begin{split}
	E=\frac{mv^2}{2} + \frac{MV^2}{2}=const\\
	\frac{mv^2}{2} + \frac{MV^2}{2}  = \frac{mv'^2}{2} + \frac{MV'^2}{2}
	\end{split}
\end{align}
\\
Закон сохранения импульса для брусков:
\begin{align}
	\label{com}
	\begin{split}
	\vec{P}=m\vec{v} + M\vec{V}=\vec{const}\\
	m\vec{v} + M\vec{V} = m\vec{v'} + M\vec{V}'
	\end{split}
\end{align}
\\
Решения уравнений \eqref{coe} и \eqref{com} для скоростей брусков после столкновения — $v'$ и $V'$ (взяты проекции скорости):
\begin{align}
	\label{sol0}
	\begin{split}
	v' = \frac{M-m}{M+m}v + \frac{2M}{M+m}V\\[6pt]
	V' = \frac{2m}{M+m}v + \frac{M-m}{M+m}V
	\end{split}
\end{align}
\\
Перепишем уравнения \eqref{sol0}б, учитывая, что $\alpha = m/M$:
\begin{align}
	\label{sol1}
	\begin{split}
	v' = \frac{\alpha-1}{1+\alpha}v + \frac{2}{1+\alpha}V\\[6pt]
	V' = \frac{2\alpha}{1+\alpha}v + \frac{1-\alpha}{1+\alpha}V
	\end{split}
\end{align}
\\
Запишем систему линейный уравнений \eqref{sol1} в матричном виде:
\begin{equation}
	\label{sol2}
	\left(
	\begin{array}{cc}
	v' \\[6pt]
	V'
	\end{array}
	\right)
	=
	\left(
	\begin{array}{cc}
	\frac{\alpha-1}{1+\alpha}&\frac{2}{1+\alpha}\\[6pt]
	\frac{2\alpha}{1+\alpha}&\frac{1-\alpha}{1+\alpha}
	\end{array}
	\right)
	\left(
	\begin{array}{cc}
	v \\[6pt]
	V
	\end{array}
	\right)
\end{equation}
\begin{equation}
	\label{matrix1}
	S=
	\left(
	\begin{array}{cc}
	\frac{\alpha-1}{1+\alpha}&\frac{2}{1+\alpha}\\[6pt]
	\frac{2\alpha}{1+\alpha}&\frac{1-\alpha}{1+\alpha}
	\end{array}
	\right)
\end{equation}

При ударе меньшего бруска о стенку вектор его скорости меняет направление(стенка в модели имеет бесконечную массу, удар абсолютно упругий). Поэтому:
\begin{align}
	\label{sol3}
	\begin{split}
	v'' = -v'\\
	V'' = V'
	\end{split}
\end{align}
Запишем систему линейных уравнений \eqref{sol3} в виде матрицы:
\begin{equation}
	\left(
	\begin{array}{cc}
	v'' \\[6pt]
	V''
	\end{array}
	\right)
	=
	\left(
	\begin{array}{cc}
	-1&0\\
	0&1
	\end{array}
	\right)
	\left(
	\begin{array}{cc}
	v' \\[6pt]
	V'
	\end{array}
	\right)
\end{equation}
\begin{equation}
	\label{matrix2}
	A=
	\left(
	\begin{array}{cc}
	-1&0\\
	0&1
	\end{array}
	\right)
\end{equation}
Таким образом, изменение скоростей брусков при столкновении между собой описывает матрица \eqref{matrix1}; удар меньшего бруска о стенку — матрица \eqref{matrix2}. Все взаимодействия в модели можно представить как последовательное применение матриц \eqref{matrix1} и \eqref{matrix2} к скоростям брусков. Это описывает произведение матриц \eqref{matrix1} и \eqref{matrix2}:
\begin{equation}
\label{matrix3}
	M=
	\left(
	\begin{array}{cc}
	\frac{1-\alpha}{1+\alpha}&\frac{-2}{1+\alpha}\\[6pt]
	\frac{2\alpha}{1+\alpha}&\frac{1-\alpha}{1+\alpha}
	\end{array}
	\right)
\end{equation}\\
Понять смысл преобразования можно произведя замены: $u=v\sqrt{\alpha}$, $u'=v'\sqrt{\alpha}$. Тогда уравнение \eqref{sol2} примет вид:
\begin{equation}
\label{sol2m}
\left(
\begin{array}{cc}
u' \\[6pt]
V'
\end{array}
\right)
=
\left(
\begin{array}{cc}
\frac{1-\alpha}{1+\alpha}&\frac{-2\sqrt{a}}{1+\alpha}\\[6pt]
\frac{2\sqrt{\alpha}}{1+\alpha}&\frac{1-\alpha}{1+\alpha}
\end{array}
\right)
\left(
\begin{array}{cc}
u \\[6pt]
V
\end{array}
\right)
\end{equation}
\begin{equation}
\label{matrix1m}
M'=
\left(
\begin{array}{cc}
\frac{1-\alpha}{1+\alpha}&\frac{-2\sqrt{a}}{1+\alpha}\\[6pt]
\frac{2\sqrt{\alpha}}{1+\alpha}&\frac{1-\alpha}{1+\alpha}
\end{array}
\right)
\end{equation}
\begin{figure}[h]
	\includegraphics[width=\linewidth]{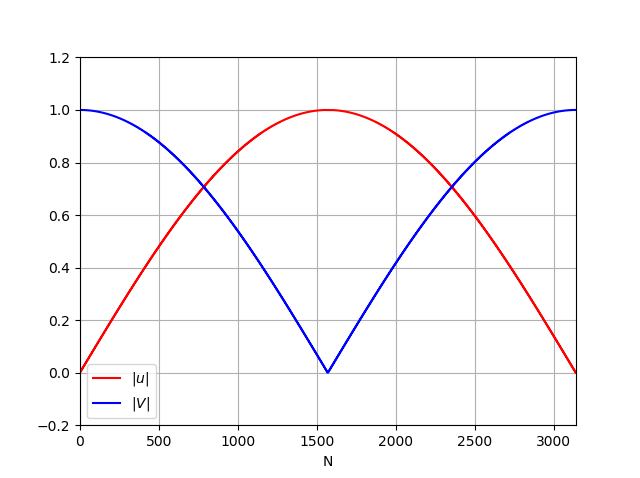}
	\caption{Модули скоростей брусков при $\alpha=10^{-6}$}
	\label{velocity}
\end{figure}
Матрица \eqref{matrix1m} является матрицей поворота вида $\left( {\begin{array}{cc}
\cos{\theta}&-\sin{\theta}\\
\sin{\theta}&\cos{\theta}\\    
\end{array}}\right)$, то есть поворачивает  вектор $(u, V)$ на угол $\theta=\arccos(\frac{1-\alpha}{1+\alpha})$. Длина вектора сохраняется, это же следует из закона сохранения энергии \eqref{coe}. Преобазуем выражение для угла $\theta$.
\begin{gather*}
\label{angle}
	\theta = \arccos(\frac{1-\alpha}{1+\alpha})\\
	\cos{\theta} = \frac{1-\alpha}{1+\alpha} = \frac{1-tg^2\frac{\theta}{2}}{1+tg^2\frac{\theta}{2}}\\
	\tg{\frac{\theta}{2}} = \sqrt{\alpha}\\
	\theta = 2\arctg{\sqrt{\alpha}}
\end{gather*}
Для определённости приняты начальные скорости $V=1$ и $v=0$ (малый брусок покоится в начальный момент времени). Столкновений больше не произойдёт при условии, что скорости брусков направлены от стены и малый брусок не догонит больший, т.е. $V<v<0$ или $V<u/\sqrt{\alpha}$. Это произойдёт, когда модуль скорости большего бруска максимально приблизиться к своей начальной скорости, но будет направлена в противоположную сторону, скорость малого шара при этом будет минимальной. Это видно по графику \eqref{velocity}, который изображает зависимости модуля скорости брусков от количества столкновений (в этом случае происходит 3141 столкновений). При столкновениях начальный вектор $(0, V_0)$ в итоге повернётся на угол, не больший $\pi$, практически развернувшись.
\begin{gather*}
	\frac{N}{2}\theta<\pi\\
	N = \Bigl\lfloor\frac{2\pi}{\theta}\Bigr\rfloor\\
	\label{xsol}
	N = \Bigl\lfloor\frac{2\pi}{2\arctg{\sqrt{\alpha}}}\Bigr\rfloor
\end{gather*}
С учетом приближени  $\arctg{\theta}\approx\theta$, можно получить решение, верное для данной модели.
\begin{equation}
\label{res}
N=\Bigl\lfloor{\pi\sqrt{\frac{M}{m}}}\Bigr\rfloor
\end{equation}
Итак, было получено решение. С учётом приближения формула \eqref{res} справедлива для всех случаев, т.к. оно выполняется при $\alpha=1$, тогда при меньших $\alpha$ приближение ещё более точно.
\section{Заключение}
Таким образом, было найдено решение задачи, впервые поставленной Галперином (2003). Столкновения рассматриваются как преобразование вектора, состоящего из скоростей брусков. Метод, предложенный Сандерсоном (2019) схож с предложенным, но в большей степени опирается на геометрические свойства.
\nocite{*}
\printbibliography
\end{document}